\documentstyle{amsppt}
\magnification\magstep1
\NoRunningHeads
\pageheight{9 truein}
\pagewidth{6.3 truein}
\baselineskip=15pt
\catcode`@=11
\def\logo@{\relax}
\catcode`@=\active
\def\af #1.{\Bbb A^{#1}}
\def\au#1.{\operatorname {Aut}\,(#1)}

\def\ring#1.{\Cal O_{#1}}
\def\pr #1.{\Bbb P^{#1}}
\def\t#1{\tilde{#1}}
\def\bq{{\Bbb Q}}
\def\pic#1.{\operatorname {Pic}\,(#1)}
\def\picq#1.{{\pic {#1}.}_{\bq}}
\def\uu{\Cal U}
\def\ol#1{\overline{#1}}
\def\ses#1.#2.#3.{0\longrightarrow #1 \longrightarrow #2 \longrightarrow #3 
\longrightarrow 0}
\def\Aut{\operatorname{Aut}}
\def\Hom{\operatorname{Hom}}

\def\mg{\overline {M}_g}

\def\ilim{\varprojlim}
\def\fmgn#1.{\overline{\Cal M}_{g,#1}}

\def\ft{{\Cal F}_{T}}
\def\mgn#1.{\overline{M}_{g,#1}}
\def\bmgn#1.{\partial{\mgn #1.}}
\def\ugn#1.{{\Cal U}_{g,#1}}
\def\sm#1.#2.{\overline{M}_{#1,#2}}
\def\sfm#1.#2.{\overline{\Cal M}_{#1,#2}}
\def\su#1.#2.{{\Cal U}_{#1,#2}}

\def\red{\operatorname{red}}
\def\ex{\Bbb E}
\def\mc#1.{\overline{NE}_1(#1)}

\define\sym#1.#2.{\operatorname {Sym}^{#1}(#2)}
\def\rpic#1.#2.#3.{\operatorname {Pic}_{#1/#2}^{#3}}

\def\sg{\operatorname {Sing}}
\def\flxy{\overline{\Cal M}^{X \cup Y}_T}
\def\lxy{\overline{M}^{X \cup Y}_T}
\def\sflxy#1.#2.{\overline{\Cal M}^{#2}_{#1}}
\def\kmg{K_1{\Cal M}_g^1}
\def\deg #1.{\operatorname {deg}(#1)}
\def\du{\coprod}
\def\deq{:=}
\def\sing{\operatorname{Sing}}
\def\supp#1.{\operatorname{supp}(#1)}
\NoBlackBoxes
\topmatter
\title Basepoint freeness for nef and big line bundles in positive
characteristic, with applications to $\fmgn n.$ and to $3$-fold MMP
\endtitle
\author Se\'an Keel
 \endauthor
\address Department of Mathematics, 
University of Texas, at Austin, Austin TX 78712 
\endaddress

\abstract
A necessary and sufficient condition is given for semi-ampleness of 
a nef and big line bundle
in positive characteristic. One application is to the geometry of the universal
stable curve over $\mg$, specifically, 
the semi-ampleness of the relative dualizing sheaf,
in positive characteristic.
An example is given which shows this, and the semi-ampleness criterion,
fail in characteristic zero. 
A second application is to Mori's program for minimal
models of $3$-folds in 
positive characteristic, namely, to 
the existence of birational extremal contractions.
\endabstract
\thanks During this research I was partially supported by NSF grant DMS-9531940
\endthanks
\endtopmatter

\heading \S 0 introduction and statement of results \endheading
A map from a variety to projective space is determined by a line bundle and
a collection of global sections with no common zeros. As all maps between 
projective varieties arise in this way, one commonly wonders
whether a given line bundle is generated by
global sections, or equivalently, if the associated linear system is basepoint free. 
Once a line bundle $L$ has a section, one expects the positive 
tensor powers $L^{\otimes n}$ to have more sections. If some such power
is globally generated, one says 
that $L$ is {\it semi-ample}.

Semi-ampleness is particularly important in
Mori's program for
the classification of 
varieties (also known as the minimal model program). Indeed a number of
the main results and conjectures --the Basepoint Free Theorem, the Abundance
Conjecture, quasi-projectivity of moduli spaces -- are explicitly issues of
semi-ampleness. I will give some details below.

There is a necessary numerical condition for
semi-ampleness. The restriction of an semi-ample line bundle to
a curve must have non-negative degree, thus: 
If the line bundle $L$ on $X$ is semi-ample, then 
$L$ is {\bf nef}, i.e. $L \cdot C \geq 0$ for every irreducible curve 
$C \subset X$. By a result of Kleiman, see \cite{Koll\'ar96,VI.2.17}, 
nefness is equivalent to the apparently stronger
condition: 
$L^{\dim Z}  \cdot Z \geq 0$ for every proper irreducible $Z \subset X$. (I
note this in relation to (0.1) below.)
Nefness does not in general imply semi-ampleness ( a
non-torsion degree zero line bundle on a curve gives a counter-example). 

The main result of this paper is a simple necessary and sufficient 
condition, in positive characteristic, 
for semi-ampleness of nef line bundles which are
close to being ample. The statement 
involves a few natural notions:

\definition{0.0 Definition-Lemma (see \cite{Koll\'ar96,VI.2.15,VI.2.16})}
A line bundle $L$ on a scheme $X$, is called
{\bf big} if 
$L^{\otimes n}$ defines a birational {\bf rational} map for
$n >> 0$.  
If $L$ is nef, and $X$ is reduced and projective over a field, 
then $L$ is big iff 
$L^{\dim X_i} \cdot X_i > 0$ for any irreducible component 
$X_i$ of $X$. 
If $L$ is semi-ample, then $L$ is big iff the associated map is
birational.
\enddefinition

Associated to a nef and line bundle is a natural locus:

\definition{0.1 Definition} Let $L$ be a nef line
bundle on a scheme $X$ proper over a field $k$. An irreducible subvariety
$Z \subset X$ is called {\bf exceptional} for $L$ if 
$L|_Z$ is not big, i.e. if $L^{\dim Z} \cdot Z =0$. If $L$ is nef
the 
{\bf exceptional locus} of $L$, denoted by $\ex(L)$, is the closure, with
reduced structure, of the union of 
all exceptional subvarieties.
\enddefinition

If $L$ is semi-ample, $\ex(L)$ is just the exceptional
locus of the associated map. It is easy to check that 
in (0.1) one does not need to take the closure, $\ex(L)$ is the union
of finitely many exceptional subvarieties, see (1.2). Of course if $L$
is not big (and $X$ is irreducible), $\ex(L) = X$. Observe that by
Nakai's criterion for ampleness, \cite{Koll\'ar96,VI.2.18}, $L$ is ample
iff $\ex(L)$ is empty.

\proclaim{0.2 Theorem} Let $L$ be a nef line bundle on a 
scheme $X$, projective over a field of positive characteristic. 
$L$ is semi-ample iff $L|_{\ex(L)}$ is semi-ample.
\endproclaim

In dimension two, (0.2), and the main ideas of its proof, 
are contained in
\cite{Artin62,2-2.11} (I received this unhappy news from Angelo Vistoli).
Similar questions were considered by Zariski, see
\cite{Zariski60}.

\proclaim{0.3 Corollary} Assumptions as in (0.2). Assume the
base field is the algebraic closure of a finite field. 
If $L|_{\ex(L)}$
is numerically trivial, then $L$ is semi-ample. In particular
this holds if $\ex(L)$ is one dimensional. If
$X$ is two dimensional, any nef line bundle which is either
big, or numerically trivial, is semi-ample. 
\endproclaim

\subhead Applications to $\fmgn n.$\endsubhead

My first application of (0.2) is to the geometry of the universal stable pointed curve
$$
\pi: \ugn n. \rightarrow \fmgn n. .
$$
Let $\Sigma \subset \ugn n.$ be the union of the $n$ universal sections. 
Using
(0.2) I prove:
\proclaim{0.4 Theorem} $\omega_{\pi}(\Sigma)$ is nef and big, and its exceptional locus
is contained in the Deligne-Mumford boundary. If the base field has positive characteristic, then
$\omega_{\pi}(\Sigma)$ is semi-ample, but this fails in characteristic zero.
\endproclaim
The nef and bigness of $\omega_{\pi}(\Sigma)$ was previously known, see e.g. 
\cite{Spiro81,pg. 56}, \cite{Viehweg89}, and \cite{Koll\'ar90,4.6}.  Note
that even nefness is not obvious: the bundle is ample on fibres by the definition of
stability, but {\it horizontal nefness} comes as a bit of a surprise. Bigness is related
to additivity of Kodaira dimension. 

The failure in characteristic zero follows from a simple example, (3.0). The same example
shows (0.2) fails in characteristic zero. It is also a counter-example to an 
interesting conjecture of Looijenga:

In Konsevich's proof of Witten's conjectures on the cohomology of $\mg$, see \cite{Kontsevich92},
a topological quotient $q: \mgn 1. \rightarrow \kmg$ plays an important role. The topological
space $\kmg$ is related to so called Ribbon Graphs, and has other applications, see
\cite{HainLooijenga96, \S 6}. In \cite{Looijenga95}, Looijenga raises the question of whether
or not $\kmg$ admits any algebraic structure. He notes that the fibres of $q$, are
algebraic, in fact (either points or) exactly the unions of exceptional curves
for $\omega_{\pi}$. This implies that 
$\omega_{\pi}$ is semi-ample iff
$\kmg$ is projective, and the quotient map
$q: \mgn 1. \rightarrow \kmg$  is the map associated to $\omega_{\pi}$. He conjectures
that these statements hold (in characteristic zero).
However
example (3.0) implies that $q$ cannot even be a morphism of schemes. See (3.6).
I don't know whether or not $\kmg$ (over ${\Bbb C}$) 
is an algebraic space (or equivalently, in terms
of definition (0.4.1) below, whether or not $\omega_{\pi}$ is EWM).

\subhead Applications to Mori's program \endsubhead

My second application of (0.2) is to Mori's program for $3$-folds, in
positive characteristic. 

Let me begin with some brief remarks on the general philosophy, so that
the statements below will be intelligible. For a detailed overview of
the program, see \cite{Koll\'ar87} and \cite{KMM87}. After I have
stated my results, I will compare them with the existing literature.

In order to build moduli spaces of varieties one looks for
a natural map to projective space, and thus for a 
natural semi-ample line bundle. On a general smooth variety, 
the only available line bundles
are $\omega_X$ and its tensor powers. If $\omega_X$ is not nef then, as noted
above, one cannot
hope for a map. Instead one looks for a birational modification which
(morally speaking) increases the nefness of $\omega_X$. 
For surfaces one blows down a $-1$ curve. One focus of Mori's program is
to generalize this procedure to higher dimensions. In this spirit, I have
the next result, which is implied by the stronger but more technical results
that follow:

\proclaim{Corollary (Existence of Extremal Contraction)} 
Let $X$ be a projective ${\Bbb Q}$-factorial normal $3$-fold defined
over the algebraic closure of a finite field.
Assume $X$ has non-negative Kodaira
dimension (i.e. $|mK_X|$ is non-empty for some
$m>0$.). If $K_X$ is not nef, then there is a surjective birational map
$f:X \rightarrow Y$ 
from $X$ to a normal projective variety $Y$, with the following
properties:
\roster
\item $-K_X$ is relatively ample.
\item $f$ has relative Picard number one. More precisely: $f$ is not an isomorphism, 
any two fibral curves are numerically equivalent, up to scaling, and
for any fibral curve $C$ the following sequence is exact:
$$
0 @>>> {\picq Y.} @>{f^*}>> {\picq X.} @>{L \rightarrow L \cdot C}>> {\Bbb Q}
@>>> 0 .
$$
\endroster
\endproclaim

(Above, and throughout the paper, a {\it fibral curve} for a map, is a curve
contained in a fibre).

A more precise statement of my results 
involves a weakening of semi-ampleness:

\definition{0.4.1 Definition} I will say that a nef line
bundle $L$ on a scheme $X$ proper over a field $k$ is
{\bf Endowed With a Map (EWM)} if there is a proper map
$f: X \rightarrow Y$ to a proper algebraic space which contracts exactly the
$L$-exceptional subvarieties -- i.e. for a subvariety $Z \subset X$,
$\dim(f(Z)) < \dim(Z)$ iff $L^{\dim Z} \cdot Z = 0$.
In particular an irreducible curve $C \subset X$ is contracted
by $f$ iff $L \cdot C = 0$. The Stein factorization of $f$ is unique, see (1.0).
\enddefinition

For indications of the relationship between EWM and semi-ample,
see (1.0).

I have a versions of (0.2-0.3) and the above corollary with EWM: whatever is stated
for semi-ampleness over a finite field, holds for EWM 
over a field of positive characteristic.
The exact statements are given in the body of the paper, see (1.9) and (1.9.1). 

\proclaim{0.5 Theorem (A Basepoint Free Theorem for big line bundles)} 
Let $X$ be a normal ${\Bbb Q}$-factorial three-fold, 
projective over a field of positive characteristic. Let $L$ be a nef
and big line bundle on $X$. 

If $L - (K_X + \Delta)$ is nef
and big for some pure boundary $\Delta$ then $L$ is EWM. If the base field
is the algebraic closure of a finite field, then $L$ is semi-ample.
\endproclaim
(A boundary is a ${\Bbb Q}$-Weil divisor $\sum a_i D_i$ with
$0 \leq a_i \leq 1$. It is called pure if $a_i < 1$ for all $i$.)

Note that when $K_X + \Delta$ has non-negative Kodaira dimension, one
does not need in (0.5) the assumption that $L$ is big, it follows from
the bigness of $L - (K_X + \Delta)$. 

Basepoint Free Theorems are related to the Existence of Extremal
Contractions, as follows: Suppose $K_X$ is not nef. Let $H$ be an
ample divisor. Let $m$ be the infimum over rational numbers 
$\lambda$ such that $K_X + \lambda H$ 
is ample. $K_X + m H$ will be nef, and should be zero on some
curve, $C$ (otherwise, at least morally, $K_X + m H$ would be ample and we
could take a smaller $m$). $L \deq K_X + m H$ is semi-ample
 by the Basepoint Free Theorem. If $f:X \rightarrow Y$ is
the associated map, then since $K_X + m H$ is pulled back from $Y$,
$-K_X$ is relatively ample. There are, however, complications. For
example, since
we take the infimum, we need to show that $m$ is rational, so that some
multiple of 
$K_X + m H$ is a line bundle. This leads to the study of the
Mori-Kleiman Cone of Curves, $\mc X.$, which is the closed convex
cone inside $N_1(X)$ (the Neron-Severi group with real coefficients)
generated by classes of irreducible curves. Specifically, one would like to
know that the edges of the cone, at least in the half space of 
$N_1(X)$ where $K_X$ is negative, are discrete, and generated
by classes of curves. In this direction I have: 

\proclaim{0.6 Proposition (A Cone Theorem for $\kappa \geq 0$)} 
Let $X$ be a normal ${\Bbb Q}$-factorial
three-fold, projective over a field. Let $\Delta$ be a boundary
on $X$.  If $K_X + \Delta$ has non-negative Kodaira dimension, then
there is a countable collection of curves $\{C_i\}$ such that
\roster
\item
$$
\mc X. = \mc X. \cap (K_X + \Delta)_{\geq 0} + \sum_i {\Bbb R} \cdot [C_i].
$$
\item All but a finite number of the $C_i$ are rational and satisfy
$0 < -(K_X + \Delta) \cdot C_i \leq 3$.
\item The collection of rays $\{{\Bbb R} \cdot [C_i]\}$ does not accumulate
in the half space $(K_X)_{< 0}$.
\endroster
\endproclaim

In characteristic zero, (0.6) is contained in \cite{Koll\'ar92,5.3}.

The proof of (0.6) is simple --since I assume that (a multiple of)
$K_X + \Delta$ is effective, the problem reduces to the cone theorem for
surfaces. 

\subhead Brief overview of related literature \endsubhead
For smooth $X$, over any base field, 
Mori's original arguments, with extensions by Koll\'ar,
give much stronger results than mine. See \cite{Mori82}, \cite{Koll\'ar91}.
The proofs are based on deformation
theory, which (at least with current technology) requires very strong assumptions
on the singularities. As one is mostly interested in smooth $X$, this may not
at first seem like a serious restriction. However even if one starts
with a smooth variety, the program may lead to singularities -- e.g. in
the above Corollary, $Y$ can have singularities, even if $X$ does not. 
Koll\'ar has been able to extend the deformation methods
to a fairly broad class of singularities, so called LCIQ singularities. 
See \cite{Koll\'ar92}. In characteristic
zero these include terminal $3$-fold singularities, the singularities
that occur in MMP beginning with smooth $X$, but this is not known 
in characteristic $p$. For other interesting applications
of Koll\'ar's main technical device, the 
Bug-Eyed cover, see \cite{KeelMcKernan95}.

In characteristic zero, with log terminal singularities, 
much stronger forms (without bigness assumptions)
of all of the above results
are known, in all dimensions. These are due to Kawamata and Shokurov,
see \cite{KMM87}. The
proofs make essential use of vanishing theorems, which fail (at least
in general) in positive characteristic. 

In one important special case, that of a semi-stable family of surfaces,
the full program is known in all (including mixed) characteristics. See
\cite{Kawamata94}. In this case the Cone and Basepoint Free theorems
are essentially surface questions, where the program is known
in all characteristics, see \cite{MiyanishiTsunoda83}.
Flips are another
(very serious) matter. 

My proof of (0.5) is based on ideas quite different from those 
of any these authors. It is a straight forward application of (0.2), and does not
use any vanishing theorems or deformation theory. Note I don't
make any singularity assumptions of the log terminal sort. 

\subhead Overview of contents \endsubhead

The proofs of (0.2-0.3) are in \S 1. Section $2$ contains technical results
about EWM and semi-ampleness used in the applications.
A counter-example to (0.2), in characteristic zero, and various
implications, are given in \S 3. I prove (0.4) in \S 4. The applications to
Mori's program are in \S 5. The section begins with the proof of (0.5). The proof
of (0.6) is in (5.5). 

{\bf Thanks:} I received help on various aspects of this paper from 
a number of people. I would like to thank in particular S. Mori, J. McKernan, 
M. Boggi, L. Looijenga, F. Voloch, Y. Kawamata, R. Heitmann, K. Matsuki, 
M. Schupsky and A. Vistoli. I would like to especially thank J. Koll\'ar for 
detailed and thoughtful (and occasionally sarcastic) 
comments on an earlier version of the paper. 

\subhead 0.9 Notation and Conventions \endsubhead
I will frequently mix the notation of line bundles and divisors. Thus
e.g. if $L$ and $M$ are line bundles, I'll write $L + M$ for $L \otimes M$.

I will often use the same symbol to denote a map, and a map 
induced from the map
by applying a functor.
I will sometimes denote the pullback $f^*(L)$ along
a map $f:X \rightarrow Y$ by $L|_X$.

$X_{\red}$ indicates the reduction of the space $X$. For a subspace
$Y \subset X$, defined by an ideal sheaf $I \subset \ring X.$, the
$k^{th}$ {\it order neighborhood} is the subscheme defined by $I^{k+1}$. 

For two Weil divisors $D,E$, I will say $D \geq E$ if the same inequality
holds for every coefficient.

All spaces considered in this paper are assumed to be separated. 

Whenever I have a base field $k$, I implicitly assume maps between
$k$-spaces are $k$-linear. The only non $k$-linear map which I consider in
the paper is the ordinary Frobenius, defined below.

{\it Frobenius Maps:}
For a scheme, $X$, of characteristic $p > 0$, and $q = p^r$, I indicate by 
$F_q: X \rightarrow X$ the (ordinary) Frobenius morphism, which is given by
the identity on topological spaces, and the $q^{th}$ power on functions. 
See \cite{Hartshorne77,IV.2.4.1}.
If $X$ is
defined over $k$, $F_q$ factors as 
$$
F_q: X \rightarrow X^{(q)} \rightarrow X
$$
where $X^{(q)} \rightarrow X$ is the pullback of $F_q$ on $\operatorname{Spec}(k)$. The
map $X \rightarrow X^{(q)}$ is called the Geometric Frobenius. It is $k$-linear. 
When $X$ is of finite type over $k$, the Frobenius, and Geometric Frobenius, are
finite universal homeomorphism. See
\cite{Koll\'ar95,\S 6}.

\heading \S 1 Proof of Main Theorem  \endheading

I begin with some simple properties of the map associated to an EWM line bundle.

\proclaim{1.0 Definition-Lemma} Let $L$ be a nef EWM line bundle on an algebraic 
space $X$, proper
over a field. I will call a proper map $X \rightarrow Y$ as in (0.4.1.) (i.e. a map 
which contracts exactly the $L$-exceptional irreducible subspaces) 
{\bf a map related to $L$}. If $f$ is such a map, and 
$f_*(\ring X.) = \ring Y.$, then $f$ is unique. I will call it the
{\bf map associated to $L$}. The associated map, $f$, has the following properties:
\roster
\item  If $f': X \rightarrow Y'$
is a proper map, which contracts any proper 
irreducible curve $C \subset X$ with $L \cdot C =0$, then $f'$ factors
uniquely through $f$. 
\item $f'$ as in (1). 
If in addition, $L$ is $f'$ numerically trivial
then the induced map $Y \rightarrow Y'$ is finite, and $f$ is the Stein factorization
of $f'$. 
\item $f$ is the Stein factorization of 
any map related to $L$.
\item Let $h:X' \rightarrow X$ be any proper map. The pullback $L|_{X'}$ is EWM, and the
associated map is the Stein factorization of $f \circ h$. 
\item $L$ is semi-ample iff $L^{\otimes m}$ is the pull back of a 
line bundle on $Y$ for some $m > 0$. 
\endroster
\endproclaim
\demo{Proof} Of course uniqueness follows from the universal property, (1), which in
turn follows from the rigidity lemma, \cite{Koll\'ar96,II.5.3}. The remaining
remarks are either easy, or contained in (1.1) and (1.3) below.
\qed \enddemo

\proclaim{1.1 Lemma} 
Let $f: X \rightarrow Y$ be a proper surjective map with geometrically connected
fibres, between algebraic spaces of finite type over a field. Assume
either that the characteristic is positive, or that 
$f_*(\ring X.) = \ring Y.$. 
If $L$ is semi-ample and $f$-numerically trivial, then, for
some $r > 0$, 
$L^{\otimes r}$ is pulled back from a line bundle on $Y$. 
\endproclaim
\demo{Proof} Let $f': X \rightarrow Y'$ be the Stein factorization of $f$.
Then by assumption, $Y' \rightarrow Y$ is a finite universal homeomorphism (the
identity in characteristic zero).
By (1.4) it is enough to show $L^{\otimes r}$ is pulled back from $Y'$.
Let $g: X \rightarrow Z$ be the map associated to $L$.
Since $L$ is pulled back from an ample line bundle on $Z$, every
fibre of $f'$ is contracted by $g$. Thus $f'$ factors through $g$ by
the rigidity lemma, \cite{Koll\'ar96,II.5.3}. \qed \enddemo

\remark{1.2} By \cite{EGAIII,4.3.4}, for proper map $f:X \rightarrow Y$,
if $f_*(\ring X.) = \ring Y.$, then $f$ has geometrically connected fibers.
\endremark

\proclaim{1.3 Corollary} Let $L$ be a nef line bundle on an algebraic
space $X$, proper
over a field. Assume $L$ is EWM and $f: X \rightarrow Y$ is the
associated map. The following are equivalent:
\roster
\item  $L$ is semi-ample.
\item $L^{\otimes r}$ is pulled back from a line bundle on $Y$,
for some $r > 0$.
\item $L^{\otimes r} $ is pulled back from an ample line bundle on $Y$,
for some $r > 0$.
\endroster
\endproclaim
\demo{Proof} (1) implies (2) by (1.1). (3) obviously implies (1).
So it is enough to show (2) implies (3). 
Assume $L = f^*(M)$ for a line bundle
$M$ on $Y$. Let  $W \subset Y$ be an irreducible subspace of dimension
$k$. Since $f$ is surjective, there is an irreducible $k$-dimensional
subspace $W' \subset X$ surjecting onto $W$. Let $d$ be the degree of
$W' \rightarrow W$. $W'$ is not $L$-exceptional, thus
$$
0 < L^k \cdot W' = d (M^k \cdot W).
$$
So $M$ is ample by Nakai's criterion, \cite{Hartshorne70}. \qed \enddemo

The next two lemmas point up the advantage of positive characteristic.

\proclaim{1.4 Lemma} Let $f:X \rightarrow Y$ be a finite universal homeomorphism
between algebraic spaces of finite type over a field $k$ of characteristic $p > 0$.
Then for some $q = p^r$ the following hold: Let $L$ be a line bundle
on $Y$.
\roster
\item For any section $\sigma$ of $f^*(L)$, $\sigma^{\otimes q}$
is in the image of 
$$
f^*: H^0(Y,L^{\otimes q}) \rightarrow H^0(X,f^*(L^{\otimes q})).
$$
\item $L$ is semi-ample iff $f^*(L)$ is semi-ample.
\item The map 
$$
f^*: \pic Y. \otimes {\Bbb Z}[1/q] \rightarrow \pic X. \otimes
{\Bbb Z}[1/q]
$$
is an isomorphism.
\item If $f^*(\sigma_1) = f^*(\sigma_2)$, for two sections
$\sigma_i \in H^0(Y,L)$, then $\sigma_1^{\otimes q} = \sigma_2^{\otimes q}$.
\endroster
\endproclaim
\demo{Proof}
By \cite{Koll\'ar95,6.6} there is a finite universal homeomorphism
$g: Y \rightarrow X$, and a $q$ as in the statement of the lemma,
such that the composition $g \circ f$ is  the Frobenius morphism, $F_{q}$ (see
(0.9)). The map induced by $F_q$ on Cartier divisors is just the $q^{th}$
power. The result follows easily. \qed \enddemo

There is also a version of (1.4) for EWM:

\proclaim{1.5 Lemma} Let $g: X \rightarrow X'$ be a finite
universal homeomorphism between algebraic spaces 
 proper over a field of positive characteristic.
A line bundle $L$ on $X'$ is EWM iff $g^*(L)$ is EWM.
\endproclaim
\demo{Proof} Let $f: X \rightarrow Z$ be the map associated to $g^*(L)$.
By \cite{Koll\'ar95,6.6} there is a pushout diagram
$$
\CD 
X  @>g>> X' \\
@V{f}VV @V{f'}VV \\
Z @>\t{g}>> Z'
\endCD
$$
with $\t{g}$ a finite universal homeomorphism. Clearly
$f'$ is a related map for $L$. \qed \enddemo

The main step in proving (0.2) is the following:

\proclaim{1.6  Proposition} Let $X$ be a projective scheme over a 
field of positive characteristic.  Let $L$ be a nef line bundle on $X$.
Suppose
$L = A + D$ where $A$ is ample and $D$ is effective and Cartier. 
$L$ is EWM iff $L|_{D_{\red}}$ is EWM. $L$ is semi-ample iff $L|_{D_{\red}}$ is
semi-ample.
\endproclaim

\demo{Proof}
Assume $L|_{D_{\red}}$ is EWM.

Let $D_k$ be the $k^{th}$ order neighborhood of $D$. By (1.5), $L|_D$ is 
EWM. Let $p:D \rightarrow Z$ be the associated map.  Let $I = I_D$.
Note that since $L|_D$ is numerically $p$-trivial,
$$
D|_{D} = L|_{D} - A|_{D}
$$
is $p$-anti-ample. Thus by Serre Vanishing (and standard exact sequences) there
exists $n > 0$ such that:

\item{I.} $R^ip_*(I^j/I^{t})=0$ for any $t \geq j \geq n$, $i > 0$. 
\item{II.} Let $J = I^s$ for some $s > 0$.
For any coherent sheaf ${\Cal F}$ on $X$, 
$R^ip_*(J^k \cdot {\Cal F} /J^{k+1} \cdot {\Cal F})$ vanishes for $k >>0$, $i>0$.
\item{III.} For $J = I^n$, $p_*(\ring ./J^t) \rightarrow p_*(\ring./J)$
is surjective for $t \geq 1$.
\smallskip

By (1.5), $L|_{D_n}$ is EWM. Let
$D_n \rightarrow Z_n$ be the associated map. By (1.0) the induced map
$p':D \rightarrow Z_n$ factors through $p$ and the induced map $Z \rightarrow Z_n$
is finite. Thus $p'$ satisfies (I-III). Replace $p$ by $p'$, and $Z$ by
the scheme-theoretic image of $p'$. I will make similar adjustments to 
$p$ later in the argument, without further remark.

By (II-III) and \cite{Artin70,3.1,6.3}, there 
is an embedding $Z_n \subset X'$ of $Z_n$ in a proper algebraic space, and
an extension of $D_n \rightarrow Z_n$ to a proper map 
$p: X \rightarrow X'$, such that $D$ is set-theoretically the
inverse image of $Z$, and such that 
$p: X \setminus D \rightarrow X' \setminus Z$ is an isomorphism.
Passing to  the Stein factorization, I may assume $p_*(\ring X.)
= \ring X'.$. It follows from (1.7) that $L$ is EWM and $p: X \rightarrow X'$
is the associated map.

Now suppose $L|_{D_{\red}}$ is semi-ample. By (1.1) 
$L^{\otimes r(k)}|_{D_k}$ is pulled back from the scheme-theoretic image
of $D_k$ in $X'$ for some 
$r(k) > 0$. Thus for some $r > 0$, $L^{\otimes r}$ is pulled back from $X'$ by
(I) and (1.10). Thus $L$ is semi-ample by (1.3). \qed \enddemo

\proclaim{1.7 Lemma} Let $L$ be a nef line bundle on
a scheme $X$ proper over a field. If $L= A + D$ with $A$ ample and
$D$ effective and Cartier, 
then $\ex(L) \subset D$. In any case $\ex(L)$ is a finite union of
exceptional subvarieties. \endproclaim
\demo{Proof} 
For the first claim, 
let $Z \subset X$ be an irreducible subvariety of
dimension $k$. If $Z \not \subset D$, then $D|_{Z}$ is 
effective and Cartier. $L|_{Z} = A|_{Z} + D|_{Z}$. Thus
$$
L^{k} \cdot Z  \geq A^k \cdot Z > 0.  
$$

I will prove the second claim by induction on the dimension of $X$. I can
assume $L$ is big.
I may write $L = A + D$ as in the statement,
by Kodaira's lemma, \cite{Koll\'ar96,VI.2.16}.
Write $D_{\red} = D_B + D_E$ where $D_B$ is the 
union of the irreducible components on which $L$ is big, and
$D_E$ is the union of the remaining components. By the first claim
$\ex(L) = D_E \cup \ex(L|_{D_B})$. \qed \enddemo

\proclaim{1.8 Lemma} Let $X$ be an algebraic space, proper over a field of
positive characteristic, and let $L$ be a nef line bundle on $X$. 
Suppose $X$ is union of closed subspaces
$X = X_1 \cup X_2$. If $L|_{X_i}$ is semi-ample (resp. EWM) for
$i = 1,2$ and $\ex(L) \subset X_1$ then $L$ is semi-ample (resp. EWM). \endproclaim
\demo{Sketch of Proof} This follows easily from \cite{Koll\'ar95,8.4} and
\cite{Artin70,6.1}. I'll state and prove a more general results in \S 2, 
which the reader may consult for a complete proof. See e.g.
(2.10.1) and (2.12). \qed \enddemo

\proclaim{1.9 Theorem} Let $L$ be a nef line bundle on a 
scheme $X$, projective
over a field of positive characteristic. 
$L$ is semi-ample (resp. EWM) iff $L|_{\ex(L)}$ is semi-ample 
(resp. EWM).
\endproclaim
\demo{Proof} I induct on the dimension of $X$. By (1.8) 
I may assume
$L$ is big. By (1.4-1.5) I may assume $X$ is reduced. 
By Kodaira's lemma, \cite{Koll\'ar96,VI.2.16}, $L = A + D$ as in (1.6). $L|_D$
is semi-ample (resp. EWM) by induction. Now apply (1.6). \qed \enddemo 

\proclaim{1.9.1 Corollary} $(L,X)$ as in (1.9). If $L|_{\ex(L)}$ is numerically
trivial (in particular if $\ex(L)$ is one dimensional) then $L$ is EWM,
and semi-ample if the base field is the algebraic closure of a finite
field. 
\endproclaim
\demo{Proof} Note any numerically trivial line bundle is EWM. (1.9.1) follows
from (1.9) and (2.16). \qed \enddemo

\proclaim{1.10 Lemma} Let $f: X \rightarrow X'$ be a proper map
between algebraic spaces, with $f_*(\ring X.) = \ring X'.$. Let
$D \subset X$ be a subspace with ideal sheaf $I \subset \ring X.$
and scheme-theoretic image $Z \subset X'$.
Let $D_k \subset X$ be the $k^{th}$ order neighborhood of $D$.
Assume 
$D$ is set-theoretically the inverse image of $Z$, and that the map
$$
f: X \setminus D \rightarrow X' \setminus Z
$$
is an isomorphism. Let $L$ be a line bundle on $X$, such that
for each $k$ there is an $r(k) > 0$ such that $L^{\otimes r(k)}|_{D_k}$ is
pulled back from (the scheme-theoretic image)  $f(D_k) \subset X'$.

If $R^1f_*(I^k/I^{k+1}) =0$ for $k >> 0$, then $L^{\otimes r}$ is
pulled back from $X'$ for some $r \geq 1$.
\endproclaim
\demo{Proof} 
Replace $L$ by a power so that $L|_D$ is pulled back from $Z$.
Choose $n$ so that 
$$
R^1f_*(I^k/I^{k+1}) = 0  \tag{1.10.1}
$$
for $k \geq n$. Let $Z_i \subset X'$ be the scheme-theoretic image of $D_i$. 
Replace $L$ by $L^{\otimes r}$, so that $L|_{D_n}$ is pulled
back from $Z_n$. I will show that 
\roster
\item $f_*(L)$ is locally free of rank one, and
\item the canonical map $f^*(f_*(L)) \rightarrow L$
is an isomorphism.
\endroster
(1-2) can be checked after a faithfully
flat extension of $X'$. The vanishing (1.10.1), and the assumption
$f_*(\ring X.) = \ring X'.$ are preserved by such  
an extension. (1-2) are local questions along $Z$. Thus I may
assume $X'$ is the spectrum of a local ring, $A$.  It follows that
$L|_{D_n}$ is trivial. By (1.6.1), I can choose, for all $i \geq 1$, 
global sections
$\sigma_i \in H^0(L \otimes D_i)$ such that 
$\sigma_i|_{D_j} = \sigma_j$ for $i \geq j$, and such that
$\sigma_1$ is nowhere vanishing. By Nakayama's Lemma, $\sigma_i$ is nowhere
vanishing and $L|_{D_i}$ is trivial, for all $i$. The collection
$\{\sigma_i\}$ induces an isomorphism
$$
\ilim H^0(L \otimes \ring D_i .) \rightarrow \ilim H^0(\ring D_i.).
$$
By the Theorem on Formal Functions, \cite{EGAIII}, the left hand side
is $H^0(L) \otimes \hat{A}$ and the right hand side is 
$H^0(\ring X.) \otimes \hat{A}= \hat{A}$, where $\hat{A}$ is the
completion of $A$ along $Z$. Hence (1).

By (1), to establish (2) I need only show surjectivity, or equivalently (in the
current local situation), that $L$ is basepoint free.
For any $k \geq n$
$$
\align
H^1(L \otimes I^{k}) \otimes \hat{A} &= \ilim_{r} H^1(L \otimes I^{k}/I^{k+r}) 
\text{ by the Theorem on Formal Functions} \\
                                   &= \ilim_{r} H^1(I^{k}/I^{k+r}) \text{ since $L|_{D_{i}}$ is 
trivial for all $i$} \\
&= 0 \text{ by (1.10.1) and induction.}
\endalign
$$
Thus $H^0(L) \rightarrow H^0(L \otimes \ring D_n.)$ is surjective. Since $L|_{D_n}$ is 
trivial, (2) follows.
\qed \enddemo

\remark{1.10.2 Remark} The proof shows that if $R^1f_*(I^k/I^{k+1}) = 0$
for all $k \geq n$, then $L^{r(n)}$ is a pull back. In characteristic $p$, by
(1.4), one need only assume $L|_D$ is pulled back, and in this case $r(n)$ can
be chosen independent of $L$. This strengthening of (1.10) was pointed out to
me by the referee. \endremark

\heading \S 2 Descending EWM or semi-ample in pushout diagrams \endheading

Suppose there
is a proper map $h: Y \rightarrow X$, and that $L$ is a nef line bundle on $X$.
If $L$ is EWM, or semi-ample, then it follows easily (see (1.0)) that the
same holds for $h^*(L)$. An important technical problem in the proofs of the main
results of the paper will be to find conditions under which 
the reverse implication holds, that is, conditions under which EWM, or semi-ampleness, 
descends from $h^*(L)$ to $L$. A simple example is (1.8). 
This section contains a number of results of this sort. 

\proclaim{2.1 Lemma} Let $i: Z \rightarrow X$ be a proper, set-theoretic surjection
between 
algebraic spaces of finite type over a perfect field $k$ of positive characteristic. Let
$f,g: X \rightarrow Y$ be maps such that $f \circ i = g \circ i$. Then there
is a finite universal homeomorphism (over $k$) $h: Y \rightarrow Y'$
such that $h \circ f = h \circ g$. 
\endproclaim
\demo{Proof} Replacing $Z$ by its scheme-theoretic image, I can assume $i$ is a closed
embedding, defined by a nilpotent ideal. For some $q = p^r$, $F_q$ (on $X$) factors
through $i$. By functorality of $F_q$, it follows that $F_q \circ f = F_q \circ g$ (here
$F_q$ is on $Y$).

$F_q$ factors as $Y \rightarrow Y^{(q)} \rightarrow Y$, where the first map is the
($k$-linear) Geometric Frobenius, and the second is the pullback of $F_q$ on $\operatorname{Spec}(k)$.
When the base field is perfect,
$Y^{(q)} \rightarrow Y$ is an isomorphism. So for $Y \rightarrow Y'$ I can take the
Geometric Frobenius. \qed \enddemo

\proclaim{2.2 Lemma} Let $L$ be a line bundle on an algebraic space
$X$ proper over a field $k$. If $k \subset k'$ is an algebraic field extension,
then $L|_{X_{k'}}$ is semi-ample (resp. EWM) iff $L$ is 
semi-ample (resp. EWM), where  $X_{k'} \deq X \times_k k'$.
\endproclaim
\demo{Proof} The semi-ample case is obvious. The EWM case follows from 
flat descent, see \cite{Artin68,7.2}. 
\qed \enddemo

{\it 2.3 Notation: 
I assume I
have a commutative diagram:
$$
\CD
{\Cal C} @>{j}>> Y \\
@V{p}VV   @V{p}VV \\
{\Cal D} @>{i}>> X \\
\endCD
\tag{2.4}
$$
which is a pushout, i.e. for any scheme $T$, the induced diagram
$$
\CD
\Hom(X,T) @>{\circ i}>> \Hom({\Cal D},T) \\
@V{\circ p}VV @V{\circ p}VV \\
\Hom(Y,T) @>{\circ j}>> \Hom({\Cal C},T)
\endCD
\tag{2.5}
$$
is a pullback diagram (of sets). I assume $i$ and $j$ are closed embeddings, $p$ is
proper and that
for any open subset $U \subset X$, the
diagram induced from (2.4) by restriction to $U$, is again a pushout. All the
spaces are algebraic spaces of finite type over a field of positive characteristic.

Let $L$ be a nef line bundle on $X$.
For any proper map $T \rightarrow X$ such that $L|_T$ is EWM, I will denote
the associated map by $g_T: T \rightarrow Z_T$, occasionally 
dropping subscripts when they are clear from context.}

\proclaim{2.6 Lemma} (Notation as in (2.3)). Assume the base field is perfect. 
Let $D \subset X$ be a reduced subspace,
and $C \subset Y$ the reduction of its inverse image. Assume ${\Cal D}$ is contained
set-theoretically in $D$. 
Assume $L|_D$ and $L_Y$ are EWM, and
that $g_Y|_C$ has geometrically connected fibres. Then $L$ is EWM. If furthermore
$L|_D$ and $L|_Y$ are semi-ample, then $L$ is semi-ample.
\endproclaim
\demo{Proof}
I will replace $L$ at various times by a positive tensor power, often without
remark. 

By (1.0), $g_C: C \rightarrow Z_C$ is the Stein factorization of
either
$$
g_Y|_C: C \hookrightarrow Y @>{g_Y}>> Z_Y
$$
or
$$
g_D \circ p|_C: C @>{p}>> D @>{g_D}>> Z_D.
$$
Let $V \deq g_Y(C) \subset Z_Y$ be the scheme-theoretic image of the
first map. Since $g_Y|_C$ has geometrically connected fibres, the induced map
$Z_C \rightarrow V$ is a finite universal homeomorphism.
Let $\ol{C},\ol{D}$ be the reductions of ${\Cal C}$ and ${\Cal D}$.

Consider first the semi-ample case. I'll use the following notation. If $T \rightarrow X$
is a map and $L|_T$ is semi-ample, then (after replacing $L$ by
a power), $L|_T$ is pulled back from an ample line bundle on $Z_T$, 
by (1.0.5). I'll write $L|_T = g_T^*(M_T)$.

(2.5), with $T = \af 1.$ gives a short exact sequence of sheaves
$$
0 \rightarrow \ring X. \rightarrow p_*( \ring Y.)\oplus \ring {\Cal D}.
@>{j^* \oplus -p^*}>>  p_*(\ring {\Cal C}.).
$$
There is an analogous exact sequence for global sections of a line bundle
$$
0 \rightarrow H^0(X,L) \rightarrow H^0(Y,L|_Y) \oplus H^0({\Cal D},L|_{\Cal D})
@>{j^* \oplus -p^*}>> H^0({\Cal C},L|_{\Cal C}).
\tag{2.7}
$$

Choose a point $x \in X$. Observe that 
$$
\align
d &\deq g_D(i^{-1}(x)) \subset Z_D, \text{ (which is $\emptyset$ if $x \not \in D$) and }\\
y &\deq g_Y(p^{-1}(x)) \subset Z_Y 
\endalign
$$
are zero-dimensional.
As $M_D$ is ample, there is a section
$\sigma^{Z_D} \in H^0(Z_D,M_D)$, non-vanishing at any point of $d$. 
Let $\sigma^D,\sigma^C$ be the pullbacks of $\sigma^{Z_D}$ to $D,C$.
By (1.4), replacing
the sections by powers, I can assume $\sigma^{C}$ 
is pulled back from a section $\sigma^V \in H^0(V,{M_Y}|_V)$. Since
$M_Y$ is ample, again replacing sections by powers, I may assume
$\sigma^V$ is the restriction of $\sigma^{Z_Y}$, non-vanishing
at any point of $y$. Let $\sigma^{Y}$ be the pullback of $\sigma^{Z_Y}$.
Let $\sigma^{\ol{D}}$ be the restrictions of
$\sigma^{D}$ to $\ol{D}$. By (1.4.1), replacing sections by powers,
I can assume $\sigma^{\ol{D}}$ is the restriction of a section
$\sigma^{\Cal D}$ in $H^0({\Cal D},L|_{\Cal D})$. 
By construction
$\sigma^{Y}$ and $p^*(\sigma^{\Cal D})$ have the same restrictions
to $\ol{C}$. Thus by (1.4.4), after replacing sections by powers, I can assume
$$
\sigma^{Y}|_{\Cal C} =p^*(\sigma^{\Cal D})
$$
and 
thus by (2.7), the sections glue to give a global section
$\sigma \in H^0(X,L)$. By construction,
$\sigma$ does not vanish at $x$. \qed 

Now for the EWM part.
Since the induced map 
$Z_C \rightarrow V$ is a finite universal homeomorphism,
by \cite{Koll\'ar95,8.4}
there is a pushout diagram
$$
\CD
Z_C @>>> V \\
@VVV   @VVV \\
Z_D @>>> Z_1
\endCD
$$
where the base is a finite universal homeomorphism. Since the left hand
column is finite, so is the right hand column.
By \cite{Artin70,6.1} there is another pushout diagram
$$
\CD
V @>>> Z_Y \\
@VVV   @VVV \\
Z_1 @>>> Z_2
\endCD
$$
where the base map is a closed embedding, and the right column is finite.
Define $f,h$ to be the compositions
$$
\align
f: Y & @>{g_Y}>>Z_Y \rightarrow Z_2 \\
h: D & @>{g_D}>> Z_D \rightarrow Z_1 \subset Z_2.
\endalign
$$
Note by construction, these are related maps for $L|_Y, L|_D$, and 
$h \circ p|_{\ol{C}} = f|_{\ol{C}}$. 
By \cite{Koll\'ar95,8.4}, there is a finite universal homeomorphism
$k:Z_2 \rightarrow Z_3$, such that $k \circ h|_{\ol{D}}$ extends to a map
$h': {\Cal D} \rightarrow Z_3$. $h' \circ p$ and $k \circ f|_{\Cal C}$ agree
on $\ol{C}$, thus by (2.1), after replacing $k$ by its composition with
another finite universal homeomorphism, I may assume
$$
p \circ h' = k \circ f \circ j \text{ (on ${\Cal C}$).}
$$
Thus by (2.5), the maps glue to give a map $g: X \rightarrow Z_3$. 
$g \circ p$ agrees with $g_Y$ up 
to composition by a finite map, and so $g \circ p$ is a related map for
for $L|_Y$. Thus by (1.14), $g$ is a related map for $L$, and $L$ is EWM.
\qed \enddemo

\proclaim{2.8 Lemma} Let $p: Y \rightarrow X$ be a proper surjection between
algebraic spaces, proper over a field $k$. Let
$g: X \rightarrow Z$ be a proper map. $L|_Y$ is EWM and $g \circ p$ is a related map, iff
$L$ is EWM and $g$ is a related map.
\endproclaim
\demo{Proof} An easy consequence of the projection formula, see e.g. the proof of (1.3).
\qed \enddemo

\proclaim{2.9 Corollary} Let $X$ be a reduced algebraic space, 
proper over a field of positive
characteristic. Suppose $X$ is union of closed subspaces $X = X_1 \cup X_2$. Let
$L$ be a nef line bundle on $X$ such that $L|_{X_i}$ is EWM for $i =1,2$. Let $g$ be
the map associated to $L|_{X_2}$. Assume $g|_{X_1 \cap X_2}$ has connected geometric
fibres. Then $L$ is EWM. If furthermore $L|_{X_i}$ is semi-ample, for $i=1,2$, then
$L$ is semi-ample.
\endproclaim
\demo{Proof} 
By (2.2) I can assume the base field is algebraically closed.
The diagram
$$
\CD
X_1 \cap X_2 @>>> X_2 \\
@VVV             @VVV \\
X_1 @>>>   X 
\endCD
$$
is a pushout diagram (where $X_1 \cap X_2$ is the scheme-theoretic intersection). So
the result follows from (2.6). \qed \enddemo

\proclaim{2.10 Lemma} Let $p: Y \rightarrow X$ be a proper surjection
between reduced algebraic spaces of finite type over a field of positive characteristic.
Let $D \subset X$ be a reduced subspace, and
$C \subset Y$ the reduction of its inverse image. Let $L$ be a line bundle on
$X$ such that $L|_D$ and $p^*(L)$ are semi-ample. 
Let $g$ be the map  associated to $p^*(L)$. Assume $g|_C$ has geometrically
connected fibres.

If $X$ is normal outside of $D$, then $L$ is semi-ample.
\endproclaim

\remark{2.10.1 Remark} Note that the connectivity assumption of (2.10) is trivially satisfied
if $p^*(L)$ is nef and big, and its exceptional locus is contained in
$C$, for in that case, any fibre of $g$ is either a single point,
or contained in $C$. \endremark

\demo{Proof of (2.10)} Let $x \in X$ be a point. 
I note first that the proof of the semi-ample
part of (2.6) yields (after replacing $L$ by a power)
sections
$$
\align
\sigma^{D} &\in H^0(D,L|_{D}) \\
\sigma^{Y}    &\in H^0(Y,p^*(L)), \text{ not vanishing at any point of $p^{-1}(x)$,} 
\endalign
$$
which have the same pullbacks to $C$.

I will construct a section of (some power of) $L$ non-vanishing at $x$. I will use 
only the existence of the above sections, and the normality of $X \setminus D$ (but no assumptions
of connectivity). 

Let $\t{X} \rightarrow X$ be the normalization of $X$, with conductors
${\Cal C} \subset {\t{X}}$, ${\Cal D} \subset X$. I may assume $D$ is the reduction
of ${\Cal D}$. 

Replacing $Y$ by a subspace, normalizing, and then taking
a normal closure (in the sense of Galois theory), 
and replacing $C$ and the sections by their pullbacks, 
I may assume $Y$ is normal, and that $p$ is generically finite and
factors as 
$$
Y @>{g}>> X' @>{j}>> \t{X} @>>> X
$$
where $X'$ is normal
$g$ is the quotient by
a finite group, and $j$ is a finite universal homeomorphism. Indeed, once
$[K(Y): K(X)]$ is a normal field extension, with Galois group $G$, 
define $X'$ to be the integral closure of $X$ in $K(Y)^G$. 

I may replace $\sigma^Y$ by the tensor power of its translates under $G$
(which replaces the restriction to $C$ by some tensor power), and so may assume
$\sigma^Y$ is $G$ invariant. Then it is pulled back from $X'$. After replacing 
sections by powers, by (1.4), I may assume $Y = \t{X}$.

By \cite{Reid94,2.1}, the diagram
$$
\CD
{\Cal C} @>>> \t{X} \\
@VVV     @V{p}VV \\
{\Cal D} @>>> X 
\endCD
\tag{2.11}
$$
is a pushout diagram (and the same is true after restricting to an open subset of $X$).
Exactly as in the proof of the semi-ample part of (2.6),
replacing sections by powers, $\sigma^{D}$ extends to $\sigma^{\Cal D}$, and
$\sigma^{\t{X}}$ and $\sigma^{\Cal D}$ have the same pullbacks to ${\Cal C}$, and 
so by (2.11) and (2.7), the sections glue to give a section of $L$, non-vanishing
at $x$. \qed \enddemo

\proclaim{2.12 Corollary} Let $X$ be a scheme projective over a field of
positive characteristic (resp. a finite field). Assume $X$ is a union of closed 
subsets $X = X_1 \cup X_2$. Let $L$ be a nef line bundle on $X$ such 
that $L|_{X_i}$ is EWM (resp. semi-ample) for $i=1,2$. Let $g_i: X_i \rightarrow Z_i$ be
the map associated to $L|_{X_i}$.  Assume that all but finitely many geometric fibres
of $g_2|_{X_1 \cap X_2}$ are connected. 
Then $L$ is EWM (resp. semi-ample).
\endproclaim
\demo{Proof} 
By (2.2) I may assume that the base field is algebraically closed.
Let $G$ be the union of the finitely many fibers of $g_2$ such that $G \cap X_1$ is 
(non-empty and) not
connected. By (2.9), 
it is enough to show that $X_1 \cup G$ is EWM (resp.
semi-ample). So I can change notation, and assume from the start that $L|_{X_2}$ is 
numerically trivial (resp. torsion by (2.16)). 
But in this case $g_i(X_1 \cap X_2)$ is zero dimensional, for $i =1,2$, 
so I can reverse the factors, repeat the argument, and reduce to the case when 
$L$ is numerically trivial. In this situation the result is obvious (resp. follows from (2.16)).
\qed \enddemo
\remark{2.12.1 Remark} Note that the connectedness condition of (2.12) holds 
in particular if $L$ is numerically trivial on $X_1 \cap X_2$, since then $g_2|_{X_1 \cap X_2}$
has only finitely many geometric fibres. 
\endremark

\proclaim{2.13 Corollary} Let $X$ be a reduced one dimensional scheme, projective
over a field of positive characteristic (resp. a finite field). Any nef line
bundle on $X$ is EWM (resp. semi-ample).
\endproclaim
\demo{Proof} 
Let $L$ be a nef line bundle on $X$.  Let $X_1$ be the union of irreducible
components on which $L$ is numerically trivial. Let $X_2$ be the union of the
remaining components. $L|_{X_2}$ is ample. $L|_{X_1}$ is EWM (resp. torsion by
(2.16)). Now apply (2.12). \qed \enddemo

\proclaim{2.14 Corollary} Notation as in (2.3). Assume the base field is
perfect, and the spaces are projective. Let $D \subset X$ be a reduced subscheme
and let $C \subset Y$ be the reduction of its inverse image.
Assume $D$ contains ${\Cal D}$ set-theoretically. Assume $L|_Y$ and
$L|_D$ are EWM (resp. semi-ample, and the base field is finite). If all but
finitely many geometric fibres of $g_Y|_C$ are geometrically connected, then $L$ is
EWM (resp. semi-ample).
\endproclaim
\demo{Proof} 
Let $G \subset Y$ be the union of the finitely fibres of $g_Y$ such that
$G \cap C$ is (non-empty and) not geometrically connected. Let $C' = C \cup G$. $L|_{C'}$
is EWM (resp. semi-ample), by (2.12.1). Clearly ${g_Y}|_{C'}$ has 
geometrically connected fibres, so (2.6) applies. \qed \enddemo

\proclaim{2.15 Corollary} Let $T$ be a reduced purely two dimensional scheme, projective
over an algebraically closed field of
positive characteristic (resp. the algebraic closure of a finite field). Let $L$ be a nef line bundle
whose restriction to each irreducible component of $T$ has numerical dimension one.
Let $p:\t{T} \rightarrow T$
be the normalization,  and let $C \subset \t{T}$ be the reduction of the conductor.
Assume $L|_{\t{T}}$ is EWM (resp. semi-ample). Assume
further that, set-theoretically, $C$ meets any generic fibre of the associated map in at most one point.
Then $L$ is EWM (resp. semi-ample).
\endproclaim
\demo{Proof} 
Immediate from (2.14) and the pushout diagram (2.11). 
\qed \enddemo

\proclaim{2.16 Lemma} Let $E$ be a projective scheme over the
algebraic closure of a finite field. Any numerically trivial
line bundle on $E$ is torsion. \endproclaim
\demo{Proof} $E$, and any fixed line bundle, are defined over a finite field $k$.
By \cite{AltmanKleiman80} the Picard functor $\rpic E. k. H.$, for
fixed Hilbert polynomial $H$, is coarsely represented by an algebraic
space, $P^{H}$,
proper, (in particular of finite type) over $k$. 
By the Riemann-Roch Theorem, \cite{Fulton84,18.3.1},
any numerically trivial line bundle
has constant Hilbert polynomial $\chi(\ring E.)$. 
Thus, since
$P^{\chi(\ring E.)}$ has only finitely many $k$-points (it is of finite type), 
the group of numerically trivial line bundles (defined over $k$) is finite. In particular
any such line bundle is torsion. \qed \enddemo

\heading \S 3 Counter-examples in characteristic zero \endheading

\subhead  Contracting the diagonal of $C \times C$ \endsubhead

Throughout this section I work over a base field $k$.

\subhead Notation for \S 3 \endsubhead
Let $C$ be a curve of genus $g$ at least $2$,
$S = C \times C$, $\pi_i$, the two projections, $\Delta \subset S$
the diagonal, and
$L = \omega_{\pi_1}(\Delta)$. Let $I = I_{\Delta}$.

\proclaim{3.0 Theorem} $L$ is nef and big. If the characteristic of 
the base field 
is positive, $L$ is semi-ample, but in characteristic
zero, $L$ is not semi-ample. (In any characteristic)
$\omega(2\Delta)$ is semi-ample, and defines a birational contraction
of $\Delta$ to a projective Gorenstein surface, with ample
canonical bundle.
\endproclaim 

\proclaim{3.0.1 Definition-Lemma}
There is a ${\Bbb Q}$-line bundle $\omega_{\pi}$
on the coarse moduli space $\ol{M}_{g,1}$, such that for any family
$f: W \rightarrow B$ of stable curves, 
and the induced map
$j: W \rightarrow \ol{M}_{g,1}$, $\omega_f$ and 
$j^*(\omega_{\pi})$
agree in $\pic W. \otimes {\Bbb Q}$. 
\endproclaim 
\demo{Proof} See \cite{Mumford77}. \qed \enddemo

\proclaim{3.1 Corollary} In characteristic zero, 
$\omega_{\pi} \in {\pic {\mg}.}\otimes {\Bbb Q}$ 
is nef and big, but 
not semi-ample, for any $g \geq 3$. \endproclaim
\demo{Proof}  
Nefness and bigness are instances of (4.4).

Let $E$ be an elliptic curve. Let $T = C \times E$.
Let $W$ be the lci surface obtained by 
gluing $S$ to $T$ by identifying $\Delta \subset S$
with the horizontal section $C \times p \subset T$. 
The first projections on each component
induces a family of stable curves of genus $g + 1$,
$f: W \rightarrow C$.
$\omega_{f}|_S = L$, see (5.3). Thus by (3.0), $\omega_{\pi}$ cannot be
semi-ample.
\qed \enddemo

\proclaim{3.2 Lemma}
$L|_{\Delta}$ is trivial, $c_1(L)^2 > 0$, 
and $L \cdot D > 0$ for any irreducible curve $D \subset S$ other
than $\Delta$. The same hold for $\omega(2\Delta)$. \endproclaim
\demo{Proof} Easy calculations using the adjunction formula. \qed \enddemo

Let $\Delta_k$ be the $k^{th}$ order neighborhood of $\Delta \subset S$.
Note that the map 
$$
\omega_C @>{\pi_1^* - \pi_2^*}>> \Omega^1_S\otimes \ring \Delta.
$$
induces an isomorphism $j: \omega_C \rightarrow I/I^2$.

\proclaim{3.3 Lemma} 
\roster
\item There is an exact sequence
$$
\ses H^1(I/I^2) . {\pic \Delta_2. }. {\pic \Delta .} .
$$
\item Let  $h: \pic C. \rightarrow H^{1,1}(C)$ be the map
$$
M \rightarrow [\pi_1^*(M) \otimes \pi_2^*(M^*)] 
\in H^1(I/I^2).
$$ 
Let
$$
g = j^{-1} \circ h: \pic C. \rightarrow H^{1,1}(C).
$$
$g(M) = c_1(M)$.
\item Let $\psi \in \Aut(S)$ switch the two factors. $\psi$
acts on $H^1(I/I^2)$ by multiplication by $-1$.
$V \in \pic \Delta_2. \otimes {\Bbb Q}$ is fixed by 
$\psi$ iff $V = \pi_1^*(L) \otimes \pi_2^*(L)$ for some
$L \in \pic C. \otimes {\Bbb Q}$. 
\endroster
\endproclaim

\demo{Proof}
(1) follows from the exact sequence of sheaves of Abelian groups
$$
1 \rightarrow I/I^2 @>{x \rightarrow 1 + x}>> {\ring \Delta_2.}^* \rightarrow
{\ring \Delta.}^* \rightarrow 1
$$
see \cite{Hartshorne77,III.4.6}.

For (2), it is enough to check the result for $M = \ring .(P)$
for a point $P \in C$. Let $P \in U \subset C$ be an affine
neighborhood of $P$, such that $P$ is cut out by $z \in \ring C.(U)$.
Let $V = C \setminus P$.  Let $\Cal U$ be the
open cover $\{U \times U,V \times V\} \cap \Delta$, of $\Delta \subset S$.
In $\pic \Delta_2.$, $\pi_1^*(\ring.(P)) \otimes \pi_2^*(\ring.(-P))$
is represented by the cocycle $z_1/z_2 \in H^1(\Cal U, \ring \Delta_2.^*)$,
where $z_i = \pi_i^*(z)$. Under the exact sequence (1) this corresponds
to the cocycle 
$$
\phi =\frac{z_1 - z_2}{z_2} \in 
H^1(\Cal U, I/I^2).
$$ 
Under the inclusion
$I/I^2 @>{d}>> \Omega^1_S \otimes \ring \Delta.$, $\phi$ maps to
the cocycle 
$$
\align
\frac{d(z_1 - z_2)}{z_2} &= d(z_1)/z_1 - d(z_2)/z_2  \\
&= (\pi_1^* - \pi_2^*)(d(z)/z) 
\endalign
$$
of $H^1(\Cal U,\Omega^1_S|_{\Delta})$ (note $z_1$ and $z_2$ are
the same in $\ring \Delta.$).
Thus $g(\ring.(P))$ is represented by the cocycle
$d(z)/z \in H^1(\Cal U,\omega_C)$, which represents $c_1(\ring .(P))$,
see e.g. \cite{Hartshorne77,III.7.4}. Hence (2). 

For (3): After tensoring with ${\Bbb Q}$, 
$\pi_1^* \otimes \pi_2^*$
induces a $\psi$ equivariant splitting of (1). By (2), $\psi$
acts on $H^1(I/I^2)$ by multiplication by $-1$. (3) follows.
\qed \enddemo

\remark{Remark} (3.3.2) is a special case of Atiyah's construction
of Chern classes (further news from A. Vistoli). See \cite{Atiyah57}. 
\endremark

\proclaim{3.4 Lemma} In characteristic zero, 
$L|_{\Delta_2}$ is non-torsion. \endproclaim
\demo{Proof} Suppose $\omega_{\pi_1}(\Delta)|_{\Delta_2}$ is 
torsion. Then obviously the same is true of $\omega_{\pi_2}(\Delta)$ and
thus 
$$
(\omega_{\pi_2}(\Delta) \otimes \omega_{\pi_1}(\Delta)^*)|_{\Delta_2} =
g(\omega_C)
$$
is torsion, contradicting (3.3.2). \qed \enddemo

\proclaim{3.5 Lemma} $\omega(2\Delta)|_{\Delta_k}$ is trivial
for all $k \geq 0$. \endproclaim
\demo{Proof} There are the usual exact sequences
$$
\ses I^k/I^{k+1} . {\ring \Delta_{k}. }. {\ring \Delta_{k-1}.}.
$$
$H^1(I^k/I^{k+1}) = H^1(\omega_C^{\otimes k})$ is trivial for
$k \geq 2$, thus by (3.2) I can assume $k =1$. $\omega(2\Delta)|_{\Delta_1}$
is in $H^1(I/I^2)$ by (3.2), and fixed by $\psi$, thus trivial
by (3.3.3). \qed \enddemo

\demo{Proof of (3.0)}
Assume (in characteristic zero) that 
$|L^{\otimes n}|$ is basepoint free for some $n > 0$. Then by
(3.2), $L^{\otimes n}$ is trivial 
on $\Delta_2$, contradicting (3.4).

By (3.5) and (3.2), as in the proof of the Example Theorem, 
$\omega(2\Delta)$ is semi-ample, and defines a contraction
$p: S \rightarrow \ol{S}$ of $\Delta$ to a 
projective normal surface $\ol{S}$. $\omega(2\Delta)$ is
the pullback of an ample line bundle, $M$ on $\ol{S}$.
${M} = K_{\ol{S}}$, since they agree away from
$p(\Delta)$. Thus $\ol{S}$ is Gorenstein, with ample dualizing sheaf.
\qed \enddemo

\proclaim{3.6 Proposition} For $g = 2h$, 
$h \geq 2$, with base field ${\Bbb C}$, 
$\kmg$ does not carry a scheme-structure so that
$f: \mgn 1. \rightarrow \kmg$ is a morphism. \endproclaim
\demo{Proof} Fix a curve $C$ of genus $h \geq 2$. Let $p: T \rightarrow C$ be
the family of stable genus $g = 2h$ curves obtained by gluing together 
two copies of $\pi_1: C \times C \rightarrow C$ along the diagonal,
$\Delta$. Let $S \subset T$ be one of the copies of $C \times C$. 

Assume $\kmg$ is a scheme, and $f$ is a (necessarily proper) morphism of schemes. 
The composition
$$
j: T @>{i}>> \mgn 1. @>{f}>> \kmg
$$
contracts $\Delta$ and is finite away
from $\Delta$. Let $q = j(\Delta)$. $q$ has an affine neighborhood, so by
taking the closure of the inverse image of a general effective principal
divisor near $q$, I can find an effective 
divisor $D \subset \mgn 1.$, which restricts to a non-trivial effective
Cartier divisor on $i(T)$, disjoint from 
$i(\Delta)$. Let $D' = i^*(D)$. $D' \subset T$ is a non-trivial effective 
Cartier divisor, disjoint from $\Delta$. 

I will follow the notation of \cite{ArbarelloCornalba87} for divisors 
on $\mgn n.$.
$i(T)$ meets only one of the irreducible boundary divisors, 
$\delta_{h,0}$, whose general element corresponds to a curve with two 
components, each of genus $h$, one of which is marked. 
$\delta_{h,0} = \pi^*(\delta_h)$, where $\pi$ is the map
$\pi: \mgn 1. \rightarrow \mg$.
Thus  
by the description of $\pic {\mgn 1.}.$, \cite{ArbarelloCornalba87,3.1}, after 
possibly replacing $D$ by a positive multiple, the class $D' \in \pic T.$ is 
$$
\omega_p^{\otimes a} \otimes p^*(N)
$$
for some integer $a$, and some line bundle $N \in \pic C.$. 
Note $\omega_p|_S = L$ (of (3.0)). Thus by (3.2), since
$D'$ is non-trivial, and disjoint
from $\Delta$, $a > 0$, and $N$ is trivial. 
But then $L$ is torsion in a Zariski
neighborhood of $\Delta$, contradicting (3.4). \qed \enddemo

\heading \S 4 The relative dualizing sheaf of the universal stable curve
\endheading
In this section I work over a base field $k$. The characteristic is
arbitrary except where noted.

\subhead Notation for \S 4 \endsubhead
I will indicate by $\fmgn n.$ the stack of $n$-pointed stable curves, and by
$\mgn n.$ its coarse moduli space. 
Let
$\pi_n: \ugn n. \rightarrow \fmgn n.$ with sections $x_1,\dots,x_n$ be 
the universal stable $n$-pointed curve. 
I will often use the same symbol for a section
of a family, and the divisor (of the total space) which is its image. I'll
also use the same symbol for a line bundle on a stack, and the associated
${\Bbb Q}$-line bundle on the coarse moduli space (see (3.0.1)). 
Let $\Sigma \subset \ugn n.$ be the union of the $n$ (disjoint) 
universal sections.
Let $L \deq L_{g,n} \deq \omega_{\pi}(\Sigma)$. 

The proof of (0.4) goes roughly as follows: First I will show that $L$
is nef and big, and 
$\ex(L) \subset \partial \mgn n+1.$. See (4.9). Then by (0.2) it is enough to
show that $L|_{\partial \mgn n+1.}$ is semi-ample. The normalizations of
the boundary components are themselves products of various $\sm g_i. n_i.$, 
and $L|_{g,n}$ restricts to {\it the analogous thing}, (4.4). Thus I can argue
by induction.

Dropping the last point $x_{n+1}$ induces the contraction functor
$\pi: \fmgn n+1. \rightarrow \fmgn n.$, which takes a stable $n+1$-pointed
curve $(C,x_1 + \dots x_{n+1})$ to the stable $n$-pointed curve
$(c(C), x_1 + \dots x_n)$. Here $c: C \rightarrow c(C)$ is the stabilization of the
$n$-pointed curve $(C,x_1 + \dots x_n)$. If the latter is stable,
$c(C) = C$. Otherwise $c: C \rightarrow c(C)$ contracts the irreducible component
(necessarily smooth and rational) of $C$ containing $x_{n+1}$ to a point. 
$\pi$ identifies
$\fmgn n+1.$ with $\ugn n.$. There is an associated diagram
$$
\CD
\ugn n+1. @>{c} >> \fmgn n+1. \times_{\fmgn n.} \ugn n. @>{p_2}>> \ugn n. \\
@V{\pi_{n+1}}VV      @V{p_1}VV                            @V{\pi_n}VV  \\
\fmgn n+1. @=    \fmgn n+1.    @>{\pi}>>       \fmgn n.
\endCD
\tag{4.0}
$$
where $c$ is the universal contraction. I will use this identification of
$\fmgn n+1.$ with $\ugn n.$ repeatedly. 
For further details see \cite{Keel92,pg. 547}.

\subhead Stratification by topological type \endsubhead
I will use the orbifold stratification of $\mgn n.$ by topological
type, which I now recall, following (with some adjustments of notation) 
\cite{HainLooijenga96,\S4}. 

Let
$B_i \subset \mgn m.$ be the locus of curves with at least $i$ singular
points. $B_i \subset \mgn m.$ has pure codimension $i$. The stratum
$B_i^0 = B_i \setminus B_{i+1}$  parameterizes 
curves with exactly $i$ singular points. 
The connected components of $B_i^0$ correspond to topological types, i.e.
equivalence classes under equisingular deformation. $B_i^0$ is an
orbifold (all I will use is that it is normal). 

There is a finite branched cover, $g:\t{B}_i \rightarrow B_i$,
a product of various $\sm i.j.$, which I will now describe. 
The map $g:\t{B}_i \rightarrow B_i$ is given roughly by normalizing a curve.
Understanding $g^*(L|_{B_i})$ is the main point in the proof of (0.4), so
I will go into some detail.

Let $T$ be a curve 
of fixed topological type, corresponding to a connected component of $B_i^0$.
Let $p: \t{T} \rightarrow T$ be the normalization of $T$. 
The labeled points $x_i$ are smooth points of $T$, and so give points
of $\t{T}$. Let $T_1,\dots,T_v$ be the irreducible components of $T$,
with normalizations $\t{T}_1,\t{T}_2, \dots \t{T}_v$.
Let 
$Y \subset \t{T}$ be $p^{-1}(\sg(T))$. Choose some ordering of the points of $Y$.
There is a fixed point free involution, $\sigma$, of $Y$ so that $T$ is recovered from
$\t{T}$ by identifying points of $Y$ according to the involution.
Let $X \subset T$ be the union of the $x_i$. Let $g_i$ be the genus of
$\t{T}_i$, and $X_i,Y_i$ the intersections of $X,Y$ with $\t{T_i}$. 
Let $\flxy$ be the product
$$
\flxy \deq \underset {1 \leq i \leq v} \to {\times} \sfm {g_i}.{X_i \cup Y_i}.
$$
and $p_i$ the projection onto the $i^{th}$ component.
Here $ \sfm {g_i}.{X_i \cup Y_i}. \deq \sfm {g_i}. r_i.$ where $r_i$ is the 
cardinality of $X_i \cup Y_i$ (this alternative notation has the advantage
of indicating how the points are distributed). 
Let
$\lxy$ be the associated coarse moduli space.
Let
$$
\uu_i \deq p_i^*(\su g_i. {X_i \cup Y_i}.)
$$
be the pullback of the universal family
$$
\pi: \su g_i. {X_i \cup Y_i}. \rightarrow \sfm {g_i}. {X_i \cup Y_i}.
$$
from the $i^{th}$ component of
$\flxy$.  
Construct a family $\ft \rightarrow \flxy$ by gluing 
together the $\uu_i$ along sections $y_j$ according to the involution
$\sigma$. There is a quotient map
$$
q: \du \uu_i \rightarrow \ft.
$$
$\ft \rightarrow \flxy$ is a family of stable $m$-pointed curves, and so induces
a map 
$g: \flxy \rightarrow \fmgn m.$ and a
commutative diagram
$$
\CD
\du \uu_i @>{q}>>  \ft @>>> \ugn m. \\
@VVV                @VVV      @V{\pi}VV \\
\flxy @=       \flxy  @>g>> \fmgn m. 
\endCD
\tag{4.1}
$$
where the right hand square is fibral.

$\ft$ has ordinary double points, 
$q: \du \uu_i \rightarrow \ft$ is the normalization, and the restriction
of the conductor to $\uu_i$ is (the union of sections) $Y_i$. Thus
by (5.3), if $h: \uu_i \rightarrow \ugn m.$ is given by the
top row of (4.1), then 
$$
h^*\omega_{\pi}(x_1 + \dots x_m) = p_i^*\omega_{\pi}(X_i + Y_i).
\tag{4.2}
$$
Let $\t{B}_i$ be the disjoint union of 
the $\lxy$, over the possible topological types $T$ (with $i$ singular points). 
$g$ induces a finite surjective map 
$g : \t{B}_i \rightarrow B_i$.
Let $\partial \lxy \subset \lxy$ be the union of the inverse images of 
the boundaries $\partial \sm {g_i}.{X_i \cup Y_i}.$ 
under the projections onto each component. Let
$\partial \t{B}_i \subset \t{B}_i$ be the union of the
$\partial \lxy$. 
Set-theoretically, 
$$
\partial \t{B}_i = g^{-1}(B_{i+1}).
\tag{4.3}
$$
The connected component of 
$B_i \setminus B_{i+1}$ corresponding to $T$ 
coarsely represents the stack obtained from
$\flxy \setminus \partial \flxy$ by forgetting the ordering on the
irreducible components. This stack is the quotient stack for the 
action of a finite group (some product of symmetric groups) on
$\flxy$, and so
in particular is smooth. Thus $B_i \setminus B_{i+1}$ is an orbifold.
All I will use is 
\proclaim{4.3.1 Lemma} $B_i \setminus B_{i+1}$ is normal.
\endproclaim

Now consider a topological type $T$ of an $n+1$-pointed
curve. Suppose $x_{n+1} \in T_j$. Let $c: T \rightarrow c(T)$ be the
contraction obtained by dropping $x_{n+1}$. 
There is a corresponding map
$$
{p_n}: \sflxy T. {X \cup Y}. \rightarrow \sflxy c(T). {X \cup Y \setminus x_{n+1}}.
$$
and a commutative diagram
$$
\CD
\sflxy T. {X \cup Y}. @>{g}>> \fmgn n+1. \\
@V{p_n}VV            @V{\pi}VV \\
\sflxy c(T). {X \cup Y \setminus x_{n+1}}. @>{g}>> {\fmgn n.} .
\endCD
$$
All but the $j^{th}$ components of $\sflxy T. {X \cup Y}.$ and
$\sflxy c(T). {X \cup Y \setminus x_{n+1}}.$ are the same. 

If $(T,x_1,\dots,x_n)$ is stable, 
then $c: T \rightarrow c(T)$ is the identity on underlying curves, and 
$p_n$ identifies $\sflxy T. {X \cup Y}.$ with
$\uu_j$, and $g: \sflxy T. {X \cup Y}. \rightarrow \fmgn n+1.$ with
$h$ of (4.2).

If $(T,x_1,\dots,x_n)$ is unstable, then the  $j^{th}$ component of 
$\sflxy T. {X \cup Y}.$ is a point, $\sfm 0. 3.$, and $p_n$ is an
isomorphism. Let $p = c(x_{n+1}) \in c(T)$. $p$ is either a singular point,
or one of the labeled points. In the first case choose some irreducible
component $p \in c(T)_s$ and a point 
$y_i$ in the normalization mapping to $p$. Then $p_n$ identifies 
$\sflxy T. {X \cup Y}.$ with the section $y_i$ of $\uu_s$. In the 
second case, say $p = x_i \in c(T)_s$. Then $p_n$ identifies
$\sflxy T. {X \cup Y}.$ with the section $x_i$ of $\uu_s$, and
$g$ with the restriction of $h$. 

\proclaim{4.4 Lemma} (Notation as above)
$$
L_{g,n}|_{\lxy} = p_j^*(L_{g_j,r_j -1})
$$
where $L_{0,2}$ indicates the trivial line bundle on
the zero dimensional space $\sfm 0. 3.$.
\endproclaim
\demo{Proof} 
This is immediate from the above identifications and (4.2), using
the adjunction formula in the unstable case to show that the left hand
side is trivial.
\qed \enddemo

\proclaim{4.5 Lemma} $L_{0,3}$ and $L_{1,1}$ are semi-ample and big.
\endproclaim
\demo{Proof} $\su 3. 0. = \pr 1.$ and $L_{0,3} = \ring {\pr 1.}.(1)$. 
The claims for $L_{1,1}$ are easily checked by considering the 
family of pointed elliptic curves given by a general pencil of plane
cubics.
\qed \enddemo

\proclaim{4.6 Lemma} Notation as in (4.0). $\omega_{\pi_{n+1}}(x_1 + \dots + x_{n})
= (p_2 \circ c)^*(\omega_{\pi_n}(x_1 + \dots + x_n))$.
\endproclaim
\demo{Proof} The left hand side is $\pi_{n+1}$ nef. In particular it is
$c$ nef.
$(\fmgn n+1. \times_{\fmgn n.} \ugn n.,x_1 + \dots + x_n)$ has
canonical singularities, by inversion of 
adjunction \cite{Koll\'ar96b,7.5} (or explicit local coordinates,
see \cite{Knudsen83}).
$c$ is birational, 
so the result follows by negativity of
contractions, see \cite{Koll\'aretal,2.19}. \qed \enddemo

\proclaim{4.7 Proposition} $L_{g,n}$ is nef. \endproclaim
\demo{Proof} First I induct on $n$ to reduce either to (4.5), or
the case $g \geq 2$ and $n =0$.  

Assume $L_{g,n}$ is nef.
By (4.6)
$$
L_{g,n+1} = (p_2 \circ c)^*(L_{g,n}) + x_{n+1} \tag{4.7.1}
$$
Since $L_{g,n+1}|_{x_{n+1}}$ is trivial by
adjunction, $L_{g,n+1}$ is nef. So I may assume $g \geq 2$ and $n=0$.

\remark{4.7.2 Remark (to be used in the proof of (4.8))} 
$$
p_2 \circ c \circ x_{n+1}: \fmgn n+1. \rightarrow \ugn n. 
\tag{4.7.3}
$$
is an isomorphism, thus 
if $d=3g-1 + n$ is the dimension of $\ugn n+1.$, then
$$
c_1(L_{g,n+1})^d \geq c_1((p_2 \circ c)^*(L_{g,n}))^{d-1} \cdot x_{n+1}
= c_1(L_{g,n})^{d-1}.
$$
Hence if $L_{g,n}$ is big, so is $L_{g,n+1}$. 
\endremark
\medskip
For (4.7) it is enough to consider a one dimensional family of stable curves. 
By (4.4) and the above reduction, I can assume $g \geq 2$, $n=0$ and the 
general fibre 
is smooth. Now the result follows from \cite{Koll\'ar90,4.6}. \qed \enddemo

\proclaim{4.8 Proposition} Let $p: {\Cal C} \rightarrow B$ with sections
$x_1,\dots,x_n$ be a family of stable $n$-pointed curves of genus $g$. Assume
$B$ is irreducible, the general fibre of $p$ is smooth, 
and the associated map $B \rightarrow \mgn n.$ is generically
finite. If $V \subset {\Cal C}$
is a subvariety surjecting onto $B$, and 
$V$ is not contained in any of the $x_i$, then $\omega_p(x_1 +\dots x_n)|_V$
is big.
\endproclaim
\demo{Proof} 
Note there are two cases. Either $V = {\Cal C}$ or $p|_V$ is generically finite.

Consider the first case. I start by inducting on $n$, to
reduce to the case of $g \geq 2$, and $n =0$, as in the proof of (4.7).
Suppose there are $n+1$ sections (and the result is known for an $n$-pointed curve).
Let $c: {\Cal C} \rightarrow {\Cal E}$ be
the contraction of $x_{n+1}$. This induces a map $B \rightarrow \mgn n.$,
which is just the composition $B \rightarrow \mgn n+1. @>{\pi} >> \mgn n.$.
After replacing $B$ by a finite cover (to deal with the fact that
$\mgn n.$ does not carry a universal family, see
\cite{Viehweg95,9.25}), I can assume there is a map $B \rightarrow B'$
so that 
${\Cal E} \rightarrow B$ is pulled back from an $n$-pointed curve
${\Cal D} \rightarrow B'$, for which the induced map $B' \rightarrow \mgn n.$
is generically finite. I have a commutative diagram
$$
\CD
{\Cal C} @>{c} >> {\Cal E} = B \times_{B'} {\Cal D} @>{p_2}>> {\Cal D} @>>> \ugn n. \\
@VVV      @V{p_1}VV                            @VVV            @V{\pi_n}VV \\
B @=    B    @>>>       B' @>>>  \fmgn n.
\endCD
$$
where the right two squares are fibral, and the left most square is 
pulled back from the left most square of (4.0). There are formulae 
analogous to (4.7.1), (4.7.3), so I can apply induction 
exactly as in (4.7.2). 

So I may assume 
$g \geq 2$ and $n =0$. 
By \cite{LaksovThorup89,4.8}
(or \cite{Viehweg77,2.10} in characteristic zero)
$$
\omega_p = p^*c_1(p_*(\omega_p)) + Z
$$
where $Z$ is an effective Cartier divisor, surjecting onto $B$ (the intersection of
$Z$ with the generic fibre is the set of Weierstrass points).
$c_1(p_*(\omega_p))$ is big by
\cite{Cornalba93,2.2}. Thus by Kodaira's lemma I may write
$c_1(p_*(\omega_p)) = A + E$, with $A$ ample, and $E$ effective and
Cartier.  Let $d-1$ be the dimension of $B$.
$\omega_p$ is nef by (4.7).
$$
c_1(\omega_p)^d \geq p^*(A)^{d-1} \cdot Z > 0.
$$
Thus $\omega_p$ is big.

Now consider the second case. After pulling back, I can
assume $V = \sigma(B)$ for a section $\sigma$, distinct from the $x_i$. I can
also assume that $B$ is normal. Thus ${\Cal C}$ is normal.

$\omega_p(\Sigma)$ is nef and big by (4.7) and the first case. Thus
by Kodaira's lemma 
$$
\omega_{p}(\Sigma)= A + E
$$ 
with $A$ ample and $E$ effective and Cartier. 

Suppose $\omega_{p}(\Sigma)|_{\sigma}$ is not big. Then by
(1.2), $\sigma$ is in the support of $E$. Thus
there exists $\lambda > 0$ so that $\lambda E = \sigma + V$, with $V$ an 
effective ${\Bbb Q}$-Weil divisor, whose support does not contain $\sigma$.
Then
$$
(1 + \lambda) \omega_{p}(\Sigma)|_\sigma = \omega_{p}(\sigma + V)|_{\sigma} + 
(\lambda A + \Sigma)|_{\sigma}
$$
the first term (on the right hand side of the above equality) is effective by
adjunction, see (5.3).
The second is big. Thus
$\omega_p(\Sigma)|_{\sigma}$ is big, a contradiction. \qed \enddemo

\proclaim{4.9 Corollary} $L_{g,n}$ is nef and big, and its exceptional
locus is contained in $\partial \mgn n+1.$.
\endproclaim

\proclaim{4.10 Theorem} If the base field has positive characteristic, then
$L$ is semi-ample.
\endproclaim
\demo{Proof}

I will proceed by induction on $m = 3g-2 + n$, the dimension of $\ugn n.$
Note $\ugn n.$ and $L$ are defined over the characteristic field, thus
(4.10) holds for $m \leq 2$ by (0.3). By (4.9) and (0.2) it 
is enough to show
$L|_{\bmgn n+1.}$ is semi-ample. 

I will prove that $L|_{B_i}$ is semi-ample by induction on $i$.
Of course $B_i$ is empty, and there is nothing to prove, for $i > m$. 
By (4.3),(4.4) and (4.9), $\ex(L|_{\t{B}_i}) \subset \partial \t{B}_i$. By
(4.9) and induction on $m$, $L|_{\t{B}_i}$ is semi-ample.
Thus $L|_{B_i}$ is semi-ample by induction on $i$, (4.3.1) and (2.10.1). \qed \enddemo

\heading \S 5 Existence of birational 
$K + \Delta$ negative extremal contractions on $3$-folds of positive
characteristic. \endheading

\subhead 5.0 Proof of (0.5) \endsubhead

By Kodaira's lemma $L = A + E$ for $A$ ample, and $E$ effective. Write
$E = N_0 + N_1 + N_2$, where $N_i$ is the sum of the irreducible components (with the same 
coefficients as in $E$) on which $L$ has numerical dimension $i$. Let $T$ be
the support of $N_1$.

The main issue will be to show that $L|_T$ is EWM (resp. semi-ample if the
base field is finite). Suppose this
has been established. Let $R_i$ be the support of $N_i$ for $i=0,2$.
Let $W = R_0 \cup \ex(L|_{R_2})$. By (1.2), 
$\ex(L) = T \cup W$. Note $\ex(L|_{R_2})$ is one dimensional, so $L|_W$ is numerically trivial
(resp. torsion by (2.16)).
Thus $L|_{\ex(L)}$ is EWM (resp. semi-ample) by (2.12.1), and so
$L$ is EWM (resp. semi-ample) by (0.2).

Now I consider $L|_T$. Let $T = \cup T_i$ be the decomposition into irreducible
components. Let $p: \t{T} \rightarrow T$ be the normalization, and
$\t{T} = \cup \t{T_i}$ the corresponding decomposition into connected components.
I will show first that $L|_{\t{T}}$ is semi-ample. Since 
$L|_{\t{T_i}}$ has numerical dimension one, it's enough to show 
$L|_{\t{T_i}}$ moves, see (5.2). This
will follow from a simple Riemann-Roch calculation. 
It will follow from adjunction that the conductor of $\t{T_i} \rightarrow T_i$ is generically
a section for the map associated to $L|_{\t{T_i}}$ (and something slightly weaker
holds for $p: \t{T} \rightarrow T$). (2.12) and (2.15) will then imply $L|_T$ is 
EWM (resp. semi-ample).

By (2.2) I may assume the base field is algebraically
closed. 

I will start off with some book keeping, to prepare for adjunction.
Let 
$M = L - (K_X + \Delta)$, which is nef and big by assumption.
Let 
$$
N = \sum_{i =1}^r a_i T_i
$$
be the decomposition into irreducible components. 
Define $\lambda_i > 0$
by $\lambda_i a_i + e_i =1$, where $e_i$ is the coefficient of $T_i$
in $\Delta$. Arrange indices so that 
$\lambda_1 \geq \lambda_2 \dots \geq \lambda_r$.
Define $\Gamma_i$ by
$$
\Delta + \lambda_i E = T_i + \sum_{j > i} T_i + \Gamma_i.
$$
By construction, $\Gamma_i$ is effective, its support does not contain
$T_i$ and 
$$
(1 + \lambda_i)L - (K_X + \sum_{j \geq i}T_j + \Gamma_i) = M + \lambda_i A 
\deq A_i
\tag{5.0.1}
$$
is ample.

Let $C_i \subset \t{T}_i$ and $D_i \subset T_i$ be the conductors.
Let $T^i = \bigcup_{j \geq i} T_j$, with normalization $\t{T^i}$ and conductor
$C^i \subset \t{T^i}$. Let $Q_i \subset \t{T_i}$ be the restriction of
$C^i$ (note this is just the restriction to a connected component). 
I will use the same symbol to indicate the integral 
Weil divisor associated to each conductor. $C_i$ is a subscheme of
$Q_i$ thus there is an inequality between Weil divisors
$Q_i \geq C_i$.

By the adjunction formula, (5.3)
$$
(K_X + \sum_{j \geq i} T_j + \Gamma_i)|_{\t{T}_i} = K_{\t{T}_i} + Q_i + R_i
\tag{5.0.2}
$$
for some effective $R_i$.
Of course (5.0.1-5.0.2) imply 
$$
K_{\t{T_i}} + Q_i + R_i = (1 + \lambda_i) L|_{\t{T_i}} - A_i|_{\t{T}_i}.
\tag{5.0.3}
$$

By the Riemann-Roch Theorem, the leading term of
$\chi(L^{\otimes r} \otimes \ring \t{T}_i.)$ is 
$$
r/2(L|_{\t{T}_i}) \cdot (L|_{\t{T}_i} - K_{\t{T}_i}).
$$
Using (5.0.3) and the fact that $L|_{\t{T_i}}$ has numerical dimension one:
$$
\align
r/2(L|_{\t{T}_i}) \cdot (L|_{\t{T}_i} - K_{\t{T}_i}) &= r/2 (L|_{\t{T_i}}) \cdot (-K_{\t{T_i}}) \\
&= r/2 (L|_{\t{T_i}}) \cdot (Q_i + R_i + A_i|_{\t{T_i}}) \\
&\geq r/2 ( L \cdot A_i \cdot T_i).
\endalign
$$
$L \cdot A_i \cdot T_i$ is strictly positive by the Hodge Index Theorem. Thus
$L|_{\t{T}}$ is semi-ample by (5.4) and (5.2). Let $f_i$ be the map associated
to $L|_{\t{T_i}}$.

By (5.0.3), $K_{\t{T_i}} + Q_i$ is negative on any generic fibre 
of $f_i$. $Q_i$, and $C_i$ have integral coefficients. Thus, by adjunction, any 
generic fibre is $\pr 1.$, and in a neighborhood of any generic fibre, either
$Q_i$ is empty, or a section. Since $Q_i \geq C_i$, the same holds for $C_i$.

Thus 
$L|_{T_i}$ is EWM (resp. semi-ample) by (2.15). Let $g_i: T_i \rightarrow Z_i$
be the associated map. 

Let 
$$
W' =\bigcup_{i=1}^{i = n-1} T_i ,
$$
with reduced structure. I will argue inductively that $L|_{W'}$ is EWM (resp. semi-ample). 
Assume this holds for $n-1$, with associated map $g: W' \rightarrow Z$, and consider
$T_n \cup W'$.
I claim that $T_n \cap W'$ meets any generic fibre of $g$ in at most one (set-theoretic)
point. The theorem follows from the claim, by (2.12).

To see the claim, suppose on the contrary, that there are points $p \neq q$ of $T_n$ along
some generic fibre, $G$.  Let $X \subset G$ be a minimal
connected union of irreducible components,  containing $p,q$ (minimal under
inclusion). Let $F$ be an irreducible component of $X$
which lies on a $T_i$ with $i$ minimal. 
By minimality of $X$, 
$$
F \cap (\{p,q\} \cup \sing(X)) 
$$
contains at least two points. Thus by the minimality of $i$, there are at least
two distinct singular points of $T^i$ along $F$, and thus $Q_i$ meets the strict
transform $\t{F} \subset \t{T}_i$ in at least two distinct points. Observe $\t{F}$ is
a generic fibre of $f_i$. But by (5.0.3), 
$$
(K_{\t{T_i}} + Q_i)|_{\t{F}} = K_{\t{F}} + Q_i|_{\t{F}} = K_{\pr 1.} + Q_i|_{\t{F}}
$$
is negative, a contradiction. \qed 

Now I turn to the lemmas used in the proof of (0.5):

\proclaim{5.2 Lemma} Let $T$ be a normal surface, projective over an algebraically
closed field. Let $L$
be a nef line bundle on $T$, of numerical dimension one. If $h^0(L^{\otimes m}) > 0$ for some
$m > 0$, then $L$ is 
semi-ample.
\endproclaim
\demo{Proof} Passing to a desingularization, I can assume $T$ is non-singular. In this
case the result is familiar, see e.g. the proof of \cite{Koll\'aretal92,11.3.1}. \qed \enddemo

\remark{5.2.1 Example} (5.2) fails if one assumes only that $T$ is integral. K. Matsuki
showed me the following counter-example: Let $C$ be a curve, and let $T$ be
obtained from $C \times C$ by gluing together two points on different fibres $F_1,F_2$
of the first projection. Take $L = \ring.(2F_1 + F_2)$ on $T$. \endremark

\proclaim{5.3 Adjunction Formula} Let $T \subset X$ be a reduced Weil
divisor on a normal variety $X$. Let $\pi: \t{T} \rightarrow T$ be
the normalization, and let $p$ be the composition
$p: \t{T} \rightarrow T \subset X$. 
Suppose $K_X + T$ is ${\Bbb Q}$-Cartier. 
Let $C \subset \t{T}$ be the Weil divisor
defined by the conductor. Then there is a canonically
defined effective $\Bbb Q$-Weil divisor $\Cal D$ on $\t{T}$, whose 
support is contained in $p^{-1}(\sing(X))$, such that
$$
K_{\t{T}} + C + {\Cal D} = p^*(K_X + T).
$$
\endproclaim
\demo{Proof} Since I am working with Weil divisors, I can remove
codimension two subsets from $T$, so may assume $T$ is Cohen-Macaulay, and
$\t{T}$ is non-singular.
By \cite{Reid94,pg. 17}, there is a canonical surjective map
$\pi^*(\omega_T) \rightarrow \omega_{\t{T}}(C)$, which is an isomorphism
wherever $T$ is Gorenstein.  
$\omega_X(T) \otimes \ring T. = \omega_T$. Thus
for each $r$ there is an induced map
$$
\omega_{\t{T}}(C)^{\otimes r} =  p^*(\omega_X(T)^{\otimes r})/\text{torsion}
\rightarrow p^*(\omega_X(T)^{[r]})/\text{torsion}
$$
which is an isomorphism wherever $X$ is non-singular. The result follows.
\qed \enddemo

Shokurov gives a formula analogous to (5.3), but without
isolating the conductor. See \cite{Shokurov91,3.1}. 
When $T$ is Gorenstein in codimension one, Corti gave 
the same formula on $T$, with the same proof. See  \cite{Koll\'aretal, 16.5}.

\remark{5.3.2 Question} I wonder if a form of (5.3) holds on $T$ itself. 
The
analog of $K$ on a reduced scheme $T$ of pure dimension $d$ is
$-\frac{1}{2} \tau_{d-1}(\ring T.)$, where $\tau$ is the Todd class, see
\cite{Fulton84,ch. 18}. In the context of (5.3) is 
$$
(K_X + T)|_T + \frac{1}{2} \tau_{d-1}(\ring T.)
$$
an effective class? This would be useful for Riemann-Roch calculations,
e.g. in the proof of (0.5). 
When $T$ is Gorenstein in codimension one, this
holds, and follows from (5.3). But it does not follow in
general from (5.3); Serre's inequality on conductor lengths
$$
n \leq 2 \delta
$$
goes in the wrong direction. See \cite{Reid94,3.2}. \endremark
\smallskip
I note one immediate corollary of (5.3): 
\proclaim{5.3.3 Corollary} Let $C \subset S$ be an integral curve on a normal
${\Bbb Q}$-factorial surface, projective over a field.
If $(K_S + C) \cdot C < 0$, then $C$ is a smooth 
curve of genus $0$.
\endproclaim
\demo{Proof} Over any singular point, the conductor has degree at least two.
Thus by degree considerations the conductor is empty and 
$\t{C} = C$ has genus $0$. \qed \enddemo

\proclaim{5.4 Lemma} Let $X$ be a projective pure dimensional scheme over a field. Let
$L$ be a nef line bundle on $X$. Then $H^{\dim(X)}(L^{\otimes r})$ is bounded
over $r > 0$. 
\endproclaim
\demo{Proof} 
$$
H^{\dim(X)}(L^{\otimes r}) = H^0(\omega_X \otimes L^{\otimes -r})^*
$$
where $\omega_X$ is the pre-dualizing sheaf, see \cite{Hartshorne77,III.7.3}. Thus
it is enough to show that for any coherent sheaf
${\Cal F}$, $H^0({\Cal F} \otimes L^{-\otimes r})$ is bounded. For this I will use
Grothendieck's D\'evissage,
\cite{Koll\'ar96,VI.2.2}, and consider the class of coherent sheaves
${\Cal F}$ with this property. Follow Koll\'ar's notation. Suppose
${\Cal F} = \ring Z.$, for an integral subscheme $Z \subset X$. 
$H^0(Z,L^{\otimes -r})$ vanishes unless $L^{\otimes r}$ is
trivial. The other two conditions of \cite{Koll\'ar96,VI.2.1} are immediate.
\qed \enddemo

\remark{Remark} For much stronger results than (5.4)  see \cite{Fujita82}. \endremark

\subhead 5.5 Cone of Curves \endsubhead

I will follow the notation of \cite{Koll\'ar96,II.4} for notions 
related to cones. In particular
$N^1(X)$ indicates the dual of the Neron-Severi group, with ${\Bbb R}$ 
coefficients. Also an extremal ray, $R$, of a closed convex cone, is a one
dimensional subcone which is extremal, i.e. if $x_1 + x_2 \in R$ then 
$x_1,x_2 \in R$. 

\definition{5.5.1 Definition} A class $h \in N^1(X)$ is
called {\bf ample } (resp. {\bf nef}) if $h$ is a strictly positive (resp.
non-negative) function on $\mc X. \setminus \{0\}$.
\enddefinition
For a line bundle on a projective variety, (5.5.1)  agrees with the usual notion
of ample, by Kleiman's criterion.

\proclaim{5.5.2 Proposition} Let $X$ be a $\Bbb Q$-factorial normal $3$-fold
projective over a field (of arbitrary characteristic). Let $\Delta$
be an ${\Bbb R}$-boundary. Let $h \in N^1(X)$ be an ample class
such that 
$$
\eta \deq K_X + \Delta + h
$$
is nef, and numerically equivalent to an effective ${\Bbb R}$-divisor, $\Gamma$.
Then the extremal subcone supported by $\eta$,
$$
\eta^{\perp} \cap \mc X.,
$$
is generated by the classes of a finite number of curves. Furthermore, 
suppose 
$$
\Gamma = \gamma + \alpha
$$
where $\alpha$ is ample, and $\gamma$ is effective. Let
$S$ be the support of $\Delta + \gamma$. Then any $\eta$ trivial
extremal ray is generated either by a curve in $\sing(S) \cup \sing(X)$,
or by a rational curve, $C$, satisfying
$$
0 < -(K_X + \Delta) \cdot C \leq 3.
$$
\endproclaim
\demo{Proof}
Let $R$ be an extremal ray with $\eta \cdot R =0$. I will argue that $R$
is of the form described. Necessarily $\gamma \cdot R < 0$. It follows
that $R$ is in the image of
$$
\mc {T}. \rightarrow \mc X.
$$
for some irreducible component $T$ of the support of $\gamma$, with $T \cdot R < 0$
(see the proof of \cite{Koll\'ar96,II.4.12}).
Let $p: \t{T} \rightarrow T$ be the minimal desingularization of the
normalization of $T$, and let $\pi: \t{T} \rightarrow X$ be the induced map.
Necessarily, $R$ is generated by the image of some $\pi^*(\eta)$ trivial
extremal ray, $J$, of $\mc \t{T}.$.

Define $a \geq 0$ so that 
$$
\Delta + a \cdot \gamma = T + E
$$
for $E$ an effective ${\Bbb R}$-divisor whose support does not contain $T$.
Note $E|_T$ has support contained in $\sing(S)$.
Let $g = h + a \alpha$. 
By (5.3) 
$$
\pi^*(K_X + T) = K_{\t{T}} + Q 
$$
for an effective class $Q$, whose image in $X$ has support contained in
$\sing(S) \cup \sing(X)$.
There is a numerical equality:
$$
(1 + a) \pi^*(\eta) = K_{\t{T}} + Q + p^*(E|_T) + \pi^*(g).
$$
Since $g$ is ample, 
$$
J \cdot (K_{\t{T}} + Q + p^*(E|_T))  < 0.
$$
Thus either $J$ is $K_{\t{T}}$ negative, or is generated by a curve mapping
into $\sing(S) \cup \sing(X)$. Suppose we are in the first case, but
not the second. By Mori's Cone Theorem (for smooth surfaces),
$J$ is generated by a smooth rational curve, $C$, with
$0 < -K_{\t{T}} \cdot C \leq 3$, and there are finitely many possible
$C$ (though the number will in general depend on $\eta$).
$T \cdot C < 0$, and $W \cdot C > 0$ for any other irreducible
component, $W$, of $S$, thus
$$
0< -(K_X + \Delta) \cdot C \leq -(K_X + T) \cdot C = 
-(K_{\t{T}} + Q) \cdot C \leq -K_{\t{T}} \cdot C \leq 3
$$
\enddemo

\remark{5.5.3 Remark: Kodaira's Lemma for ${\Bbb R}$-classes} 
If one could show that a nef class $\eta$ 
with $\eta^{\dim X} > 0$ satisfied Kodaira's lemma, i.e. has
an expression $h + E$ with $h$ ample, and $E$ effective, 
then (5.5.2) would imply that 
the supported extremal subcone is finite rational
polyhedral. In
dimension two this is easy, see \cite{Koll\'ar96,4.12}. It is also known,
in all dimensions, in characteristic zero; the proof, due to
Shokurov, is an application of Kawamata-Viehweg Vanishing, see
\cite{Shokurov96}. This generalized
Kodaira lemma has some other interesting implications, see
\cite{KeelMcKernan96,\S 2}.
\endremark

\demo{Proof of (0.6)}
Let $K_X + \Delta = \gamma$, for an effective class
$\gamma$. Let $R$ be an extremal ray, with $(K_X + \Delta) \cdot R < 0$.
Let $\eta \in N^1(X)$ be a nef class supporting $R$. Then by compactness
of a slice of $\mc X.$, \cite{Koll\'ar96,II.4.8}, after possibly 
replacing $\eta$ by a positive multiple, 
$$
h \deq \eta - (K_X + \Delta)
$$
is ample. (1-3) now follow from (5.5.2), see e.g. the proof of
\cite{Koll\'ar96,III.1.2}.
\qed \enddemo

\remark{5.5.4 Remark} It is natural to expect, under some singularity assumptions,
that the extremal rays in (0.6) are all generated by smooth rational
curves. In characteristic zero this follows from Kawamata-Viehweg
Vanishing applied to the corresponding extremal contraction. In characteristic
$p$, the proof of (0.5), plus some easy analysis in the case when the
exceptional locus is a surface contracted to a point, shows that the ray is generated by a 
(possibly singular) rational curve, except possibly in the case of a small
contraction. For a small contraction, 
if the singularities are isolated $LCIQ$ and the ray is
$K_X$ negative, then any generating curve is rational by \cite{Koll\'ar92,6.3}.
\endremark

\bigskip
Mathematics Department,
University of Texas at Austin
Austin, Texas 78712 USA

email: keel\@math.utexas.edu

\end